\documentstyle[12pt]{article}
\begin{document}
\author{S.V. Ludkovsky and B. Diarra.}
\title{Profinite and Finite Groups Associated with
Loop and Diffeomorphism Groups of Non-Archimedean Manifolds.
}
\date{08 December 2000}
\maketitle
\section{Introduction.} 
The importance of such groups in the non-Archimedean functional 
analysis, representation theory and mathematical physics is clear and
also can be found in the references given below
\cite{aref,luseabm,lutmf,luladg,mensk,roo,roosch}.
This article is devoted to one aspect of such groups:
their structure from the point of view of the $p$-adic 
compactficiation (see also about 
Banaschewski compactification in \cite{roo}).
This also opens new possibilites for studying their 
representations as restrictions of representations of 
$p$-adic compactifications, which are constructed below 
such that they also are groups.
\par At first we remind basic facts and notations,
which are given in detail in references \cite{luseabm,lutmf,lubp,roo,sch1}.
For a diffeomorphism group $Diff(M)$ of a Banach manifold
over a local field $\bf K$ there are clopen subgroups $W$ such that
they contain a sequence of profinite subgroups $G_n$ with
$G_n\subset G_{n+1}$ for each $n\in \bf N$ and $\bigcup_nG_n$
is dense in $W$. A loop group $L_t(M,N)$ is defined as a quotient space
of a family of mappings $f: M\to N$ of class $C^t$ of one Banach manifold  $M$ into another $N$ over the same local field $\bf K$
such that $\lim_{z\to s}(\bar \Phi ^vf)(z;h_1,...,h_n;\zeta _1,...,
\zeta _n)=0$
for each $0\le v\le t$, where $M$ and $N$ are embedded into the corresponding
Banach spaces $X$ and $Y$, $cl(M)=M\cup \{ s \} $, $cl(M)$ and $N$
are clopen in $X$ and $Y$ respectively, $0\in N$, 
$(\bar \Phi ^vf)(z;h_1,...,h_n;\zeta _1,..,\zeta _n)$ 
are continuous extensions of difference quotients, $z\in M$, 
$h_1,...,h_n$ are nonzero vectors in $X$, 
$\zeta _1,...,\zeta _n\in \bf K$ such that 
$z+\zeta _1h_1+...+\zeta _nh_n\in M$, $n=[v]+sign \{ v \} .$ 
$p$-adic completions of clopen subgroups $W$ of loop groups $G$
and diffeomorphism groups $G$ are considered.
In the case of the diffeomorphism group the $p$-adic completion
produces weakened topology on $W$ relative to which it remains
a topological group. In the case of the loop group the $p$-adic
completion produces a new topological group $V$ in which the initial group
$W$ is embedded as a dense subgroup such that $V\ne W$.
The topology on $W$ inherited from $V$ is weaker than the initial one.
For the compact manifold $M$ in the case of the diffeomorphism group
the $p$-adic completion of $W$ produces profinite group.
For the locally compact manifolds $M$ and $N$ in the case of the loop group
$L^t(M,N)$ the $p$-adic completion of $W$ produces its embedding
into ${\bf Q_p}^{\bf N}$. When $W$ is bounded relative
to the corresponding metric in $L^t(M,N)$, then $W$ is embedded
into ${\bf Z_p}^{\bf N}$. The group $Diff(M)$ is perfect and simple,
on the other hand, the group $L^t(M,N)$ is commutative.
The notation given below and the corresponding definitions are given 
in detail in \cite{luseabm,lubp}.
\section{$p$-Adic Completion of Diffeomorphism Groups.}
\par {\bf 2.1. Notations and Remarks.} Let $N$ be a compact manifold over 
a local field $\bf K$, that is, $\bf K$ is a finite algebraic extension
of the field of $p$-adic numbers $\bf Q_p$ \cite{wei}.
Let also $N$ be embedded into $B({\bf K^{\sf n}},0,1)$ 
as a clopen (closed and open at the same time) subset \cite{boum,luum985}, 
where ${\sf n}\in \bf N$, $B(X,y,r):= \{ z: z\in X; d_X(y,z)\le r \} $
denotes a clopen ball in a space $X$ with an ultramentric $d_X$. 
The ball $B({\bf K^{\sf n}},0,1)$ has the ring structure with coordinatewise
addition and multiplication, in particular $B({\bf Q_p},0,1)=\bf Z_p$
is the ring of entire $p$-adic numbers.  
The ring $B({\bf K^{\sf n}},0,1)$ is homeomorphic with the projective limit
$B({\bf K^{\sf n}},0,1)=pr-\lim_k{\bf S_{p^k}}^{\sf n}$, 
$\bf S_{p^k}$ is a finite ring
consisting of $p^{kc}$ elements such that $\bf S_{p^k}$ is equal
to the quotient ring $B({\bf K},0,1)/B({\bf K},0,p^{-k})$, 
${\bf S_{p^k}}^{\sf n}$ is a product of $\sf n$ copies of ${\bf S_{p^k}}$,
$c$ is a dimension $dim_{\bf Q_p}{\bf K}$ of $\bf K$ over $\bf Q_p$
as a $\bf Q_p$-linear space.
In particular $B({\bf Q_p},0,1)/B({\bf Q_p},0,p^{-k})
={\bf Z_p}/p^k{\bf Z_p}=\bf F_{p^k}$         
is a finite ring consisting of $p^k$ elements,
$\bf Z_p$ is the ring of $p$-adic integer numbers,
$aB:= \{ x: x=ab, b\in B \} $ for a multiplicative group $B$ 
and its element $a\in B$, $k\in \bf N$ \cite{roo,wei}. 
For each $m\ge k$ there are the following quotient mappings
(ring homomorphisms): $\pi _m: B({\bf K},0,1)\to {\bf S_{p^m}}$ and 
$\pi ^m_k: {\bf S_{p^m}}\to \bf S_{p^k}$. This induces the quotient 
mappings $\pi _m: N\to N_m$ and $\pi ^m_k: N_m\to N_k$, where
$N_m\subset {\bf S_{p^m}}$, $\pi ^m_k\circ \pi _m=\pi _k$.
\par Let now $M$ and $N$ be two analytic compact manifolds embedded
into $B({\bf K^{\sf m}},0,1)$ and $B({\bf K^{\sf n}},0,1)$ respectively
as clopen subsets and $f\in C^t(M,N)$, where $C^t(M,N)$ denotes 
the space of functions $f: M\to N$ of class $C^t$, $t\ge 0$.
For an integer $t$ it is a space of $t$-times continuously 
differetiable functions in the sence of partial difference quotients
(see in detail: \cite{luseabm,lubp,sch1}).
Then $f=pr-\lim_kf_k$, where $f_k:=\pi _k\circ f$.
We inroduce the following notation: 
$C^t(M,N_k):=\pi _k \circ C^t(M,N) = \{ f_k: f\in C^t(M,N) \} $,
hence $C^t(M,N)=pr-\lim_kC^t(M,N_k)$ algebraically without taking 
into account topologies (or the limit of the inverse 
sequence, see \S 2.5 \cite{eng} and \S \S 3.3, 12.202 \cite{nari}). 
Each function $f\in C^t(M,N)$ has a 
$C^t(B({\bf K^{\sf m}},0,1),{\bf K^{\sf n}})$-extension 
by zero on $B({\bf K^{\sf m}},0,1)$, hence it has the following decomposition 
$f=\sum_{l,m}f^l_m{\bar Q}_me_l$
in the Amice's basis ${\bar Q}_m$, where $e_l$ is the standard 
orthonormal basis in $\bf K^{\sf n}$ such that $e_l=(0,...,0,1,0,...)$
with $1$ in the $l$-th place, $m\in \bf Z^{\sf n}$,
$m_l\ge 0$, $m=(m_1,...,m_{\sf n})$, $f^l_m\in \bf K$
are expansion coefficients such that $\lim_{l+|m|\to \infty }
|f^l_m|_{\bf K}J(t,m)=0$, ${\bar Q}_m$ are polynomials
on $B({\bf K^{\sf m}},0,1)$ with values in $\bf K$.
The space $C^t(M,N)$ is supplied with the uniformity
inherited from the Banach space $C^t({\bf K^{\sf m}},{\bf K^{\sf n}})$,
$J(t,m):=\| {\bar Q}_m \|_{C^t({\bf K^{\sf m}},{\bf K})}$. 
\par {\bf 2.2. Lemma.} {\it Each $f\in C^t(M,N)$
is a projective limit $f=pr-\lim_kf_k$ of polynomials
$f_k=\sum_{l,m}f^l_{m,k}{\bar Q}_{m,k}e_l$ on rings
$\bf S_{p^k}^{\sf m}$ with values in $\bf S_{p^k}^{\sf n}$,
where $f^l_{m,k}\in \bf S_{p^k}$ and ${\bar Q}_{m,k}$
are polynomials on $\bf S_{p^k}^{\sf m}$ with values in
$\bf S_{p^k}$.} 
\par {\bf Proof.} In view of \S 2.1 $f_k=\pi _k\circ f$
and $\pi _k\circ f(x)=\sum_{l,m}(\pi _k (f^l_m)(\pi _k
{\bar Q}_m(x))e_l$, since $\pi _k$ is a ring homomorphism
and $\pi _k(e_l)=e_l$.
Then $\pi _k(ax^m)=a_kx^m(k)$ for each $a\in \bf K$ and
$x\in B({\bf K^{\sf m}},0,1)$, where 
$x^m=x_1^{m_1}...x_{\sf m}^{m_{\sf m}},$ $a_k=\pi _k(a)$
with $a_k\in \bf S_{p^k}$ and $x^m(k)=\pi _k(x^m)$
with $x(k)\in \bf S_{p^k}^{\sf m},$ consequently,
$\pi _k({\bar Q}_m(x))={\bar Q}_{m,k}(x(k)).$
The series for $f_k$ is finite, since $\pi _k(a)=0$ for each 
$a\in \bf K$ with $|a|<p^{-k}$ and  $\lim_{l+|m|\to \infty }
|f^l_m|_{\bf K}J(t,m)=0$.
\par {\bf 2.3. Corolary.} {\it The space $C^t(M,N_k)$ is isomorphic with
the space $N_k^{M_k}$ of all mappings from $M_k$ into $N_k$.
Moreover, $({\bf S_{p^k}^{\sf n}})^{({\bf S_{p^k}^{\sf m}})}$
is a finite-dimensional space over the ring $\bf S_{p^k}$.}
\par {\bf Proof.} From \S 2.2 it follows that there is only a finite number
of $\bf S_{p^k}$-linearly independent polynomilas
${\bar Q}_{m,k}(x(k))$ (that is, in the module of the ring
$\bf S_{p^k}$), since the rings ${\bf S_{p^k}^{\sf m}}$
and $\bf S_{p^k}$ are finite, also $z^a=z^b$
for each natural numbers $a$ and $b$ such that 
$a=b\mbox{ } (mod\mbox{ }(p^k))$ and each $z\in 
\bf S_{p^k}$. The space $C^t(M,N_k)$
is isomorphic with $N_k^{M_k}$, since $M_k$ and $N_k$ are discrete.
\par {\bf 2.4. Corollary.} {\it $\pi _k \circ Diff^t(M)$
is isomorphic with the symmetric group $S_{\sf n_k}$,
where $\sf n_k$ is the cardinality of $M_k$.}
\par {\bf Proof.} If $h\in Diff^t(M)$, then $h_k(M_k)=M_k$, since
$h(M)=M$. In view of \S 2.3 $\pi _k \circ Diff^t(M)$
is isomorphic with the following group $Hom(M_k)$ of all 
homeomorphisms $h_k$ of $M_k$, that is, bijective surjective mappings
$h_k: M_k\to M_k$. Using an enumeration of elements of $M_k$ we get
an isomorphism of $Hom(M_k)$ with $S_{\sf n_k}$.
\par {\bf 2.5.} Let $C_w(M,N):=pr-\lim_kN_k^{M_k}$
be a uniform space of continuous mappings $f: M\to N$ supplied with
a uniformity inherited from products of uniform spaces
$\prod_{k=1}^{\infty }N_k^{M_k}$ (see also \S 8.2 \cite{eng}). 
The spaces $C^t(M,N)$ and $C_w(M,N)$ are subsets of
$\bf K$-linear spaces $C^t(M,{\bf K^{\sf n}})$ and
$C^0(M,{\bf K^{\sf n}})$ respectively. We consider
algebraic structures of subsets of the latter $\bf K$-linear
spaces as inherited from them.
\par {\bf Corollary.} {\it The space $C^t(M,N)$ is not algebraically 
isomorphic with $C_w(M,N)$, when $t>0$. The uniform space $C_w(M,N)$ is
uniformly isomorphic with $C^0(M,N)$, when the latter space is supplied 
with a weak uniformity inherited from $C^0(M,{\bf K^{\sf n}})$.
The space $C_w(M,N)$ is compact.}
\par {\bf Proof.} In view of \S 2.5 \cite{eng} we have that
$C^0(M,N)$ and $C_w(M,N)$ coincide algebraically, since the connecting 
mappings $\pi ^m_n$ are uniformly continuous for each $m\ge n$.
The space $C^0(M,{\bf K^{\sf n}})$ is $\bf K$-linear and its uniformity
is completely defined by a neighbourhood base of zero.
The set of all evaluation mappings in points of $M$ produces the
base of the topology of $C^0(M,{\bf K^{\sf n}})$. In its weak topology
the latter space is isomorphic with the product 
$\prod_{x\in M}{\bf K^{\sf n}}={\bf K}^{card (M)}$, since
$card(M)=card ({\bf R})=\sf c$. Then $C^0(M,N)$ and $C_w(M,N)$
have embeddings into $B({\bf K},0,1)^{card (M)}$ as closed bounded 
subspaces. The latter space is uniformly homeomorphic with
$pr-\lim_k({\bf S_{p^k}})^{M_k}$, which is compact by 
the Tychonoff theorem 3.2.4 \cite{eng}. 
Since $C^0(M,N)\ne C^t(M,N)$ for $t>0$,
then $C_w(M,N)$ and $C^t(M,N)$ are different algebraically.
\par {\bf 2.6.} Let $Diff_w(M):=pr-\lim_k Hom(M_k)$
be supplied with the uniformity inherited from $C_w(M,M)$.
The group $Diff_w(M)$ is called the $p$-adic compactification 
of $Diff^t(M)$. The following theorem shows that this 
terminology is justified.
\par {\bf Theorem.} {\it $Diff_w(M)$ is the compact topological group
and it is the compactification of $Diff^t(M)$ in the weak topology.
If $t>0$, then $Diff^t(M)$ does not coincide with $Diff_w(M)$.}
\par {\bf Proof.} Since $Diff^t(M)\subset C^t(M,M)$, then 
$Diff^t(M)$ has the corresponding embedding into $C_w(M,M)$.
Since $C_w(M,M)$ is compact and $Hom(M)$ is a closed subset in $C_w(M,M)$, 
then due to Corollary 2.5 $Hom(M)\cap C_w(M,M)=Diff_w(M)$ is compact.
The space $C^t(M,M)$ is dense in $C^0(M,M)$, consequently, 
$Diff^t(M)$ is dense in $Diff_w(M)$. If $t>0$, then $Diff^t(M)\ne Hom(M)$,
hence $Diff^t(M)$ and $Diff_w(M)$ do not coincide algebraically.
It remains to verify, that $Diff_w(M)$ is the topological group
in its weak topology. If $f, g \in C^t(M,N)$, then
$\pi _k({\bar Q}_m(g(x)))={\bar Q}_{m,k}(g_k(x(k))$, consequently,
$\pi _k(f\circ g)=\sum_{l,m}\pi _k(f^l_m){\bar Q}_{m,k}(g_k(x(k))e_l$
and inevitably $(f\circ g)_k=f_k\circ g_k.$ 
On the other hand $\pi _k(x)=x(k)$, hence $\pi _k(id(x))=id_k(x(k)),$
where $id(x)=x$ for each $x\in M$. Therefore, for $f=g^{-1}$ we have
$(f\circ g)_k=f_k\circ g_k=id_k$, hence $\pi _k(g^{-1})=g_k^{-1}$.
The associativity of the composition 
$(f_k\circ g_k)\circ h_k=f_k\circ (g_k\circ h_k)$ of all functions
$f_k, g_k, h_k\in Hom(M_k)$ together with others properties given above
means, that $Diff_w(M)$ is the algebraic group, since $f=pr-\lim_kf_k$,
$g=pr-\lim_kg_k$ and $h=pr-\lim_kh_k$ also satisfy the associativity
axiom, each $f$ has the inverse element $f^{-1}(f(x))=id$ and $e=id$ 
is the unit element. By the definition of the weak 
topology in $Diff_w(M)$ for each neighbourhood of  
$e=id$ in  $Diff_w(M)$ there exists $k\in \bf N$ and a subset
$W_k\subset Hom(M_k)$ such that $e_k\in \pi _k^{-1}(W_k)\subset W$.
But $Hom(M_k)$ is discrete, hence there are $e_k\in V_k\subset Hom(M_k)$ 
and $e_k\in U_k \subset Hom(M_k)$ such that $V_kU_k\subset W_k$,
hence there are neighbourhoods $e\in V\subset Diff_w(M)$
and $e\in  U\subset Diff_w(M)$ such that $VU\subset W$, where
$V=\pi _k^{-1}(V)$, $U=\pi _k^{-1}(U)$ and $VU= \{ h: h=f\circ g,
f\in V, g\in U \} $. If $W'$ is a neighbourhood of $f^{-1}$, 
then $V:=W'f^{-1}$ is the neighbourhood of $e$ and there exists 
$k\in \bf N$ such that $\pi _k^{-1}(e_k)=:U\subset V^{-1}$,
since $e_k^{-1}=e_k$ and $\pi _k$ is the homomorphism.
Therefore, $fU:=W$ is the neighbourhood of $f$ such that
$W^{-1}\subset W'$, that demonstrates the continuity of 
the inversion operation $f\mapsto f^{-1}$.
\par {\bf 2.7. Notes.} Each projection $\pi _k: C^t(M,{\bf K^{\sf n}})
\to ({\bf K_k^{\sf n}})^{M_k}$ produces the quotient metric $\rho _k$
in the $\bf K_k$-module $({\bf K_k^{\sf n}})^{M_k}$
such that $\rho _k(f_k,g_k):=inf_{z, \pi _k(z)=0} \| f-g+z \| _{C^t
(M,{\bf K^{\sf n}})},$
where ${\bf K_k}:={\bf K}/B({\bf K},0,p^{-k})$ is the quotient ring
and $\pi _k$ is induced by such quotient mapping 
from $\bf K$ onto $\bf K_k$.
If $C^t(M,{\bf K^{\sf n}})$ embed into $\prod_k
\pi _k(C^t(M,{\bf K^{\sf n}}))$ and supply the latter space with
the box topology given by the following norm
$\| f-g \| ':=\sup_k\rho _k(f_k,g_k)$, then it produces the uniformity in
$C^t(M,{\bf K^{\sf n}})$ equivalent with the initial one.
\par Theorem 2.6 means that the $p$-adic completion $Diff_w(M)$
is the profinite group, that is, it is projective limit of 
finite groups $Hom(M_k)$. If the compact manifold $M$ is 
decomposed into the disjoint union $M=\bigcup_iB({\bf K^{\sf m}},x_i,r_i)$ 
of clopen balls, then orders of the latter groups are divisible by $(p^a)!$, 
where $a=\sum_il_i$, $l_i=k-\max_l \{ l: p^l
\le r_i \} $, $x_i\in B({\bf K^{\sf m}},0,1)$, $0<r_i\le 1$,
since $card(M_k)$ is divisible by $p^a$.
Then the representations of symmetric groups known from the classical works
of A. Young and H. Weyl \cite{litl,weyl}
with the help of the projective limit decompositions
produce finite-dimensional representations of the diffeomorphism groups.
\section{$p$-Adic Completion of Loop Groups.}
\par {\bf 3.1.} Let as in \S 2.1 $\bar M$ and $N$ be two compact manifolds
and $Diff^t_0({\bar M})$ be a subgroup in $Diff^t({\bar M})$ of all elements 
$\psi \in Diff^t({\bar M})$
such that $\psi (s_0)=s_0$, where $s_0$ is a marked point in $\bar M$.
We denote shortly by $C^t_0(M,N)$ a subspace in $C^t({\bar M},N)$
of all elements $f \in C^t({\bar M},N)$ such that
$\lim_{|\zeta _1|+...+|\zeta _n|\to 0}{\bar \Phi }^v(f-w_0)
(s_0;h_1,...,h_n;\zeta _1,...,\zeta _n)=0$ for each $v \in \{ 0,1,...,
[t],t \} $, $n=[v]+sign \{ v \} $, where $M={\bar M}\setminus s_0$ 
and $w_0({\bar M})= \{ y_0 \} $
(see \S 2.6 \cite{lubp}). 
\par {\bf Theorem.} {\it Let $\Omega _{\xi }(M,N)$ be commutative
loop monoids, then the quotient mappings $\pi _k$ induce
the corresponding inverse sequence  $\{ \Omega (M_k,N_k): k\in {\bf N} \} $
such that $\Omega ^w(M,N):=pr-\lim_k \Omega (M_k,N_k)$ 
is the commutative compact topological monoid, where
$\pi _k: \Omega _{\xi }(M,N)\to \Omega (M_k,N_k)$, 
$\pi ^l_k: \Omega (M_l,N_l)\to \Omega (M_k,N_k)$
are surjective mappings for each $l\ge k$,
$\Omega (M_k,N_k)=\{ f_k: f_k\in N_k^{M_k}, f_k(s_{0,k})=y_{0,k} \} 
/K_{\xi ,k}$, $K_{\xi ,k}$ is an equivalence relation induced by
an equivalence relation $K_{\xi }$. Moreover, $\Omega ^w(M,N)$ is
a compactification of $\Omega _{\xi }(M,N)$.}
\par {\bf Proof.} In view of Corollary 2.3
$\pi _k(C^{\xi }_0(M,N))$ is isomorphic with 
$\{ f_k: f_k\in N_k^{M_k}, f_k(s_{0,k})=y_{0,k} \} $,
where the quotient mapping is denoted by $\pi _k$ both for $M$ and $N$,
since it is induced by the same ring homomorphism $\pi _k: {\bf K}
\to {\bf K}/B({\bf K},0,1)$, $s_{0,k}:=\pi _k(s_0)$
and $y_{0,k}:=\pi _k(y_0)$. Then $\pi _k(Diff^t_0(M))$
is isomorphic with $Hom_0(M_k):=\{ \psi _k: \psi _k\in Hom(M_k),
\psi _k(s_{0,k})=s_{0,k} \} $. 
All of this is also applicable 
with the corresponding changes to classes of smoothness
$C^{\xi }$ (or $C({\xi })$ in the notation of \cite{lubp}), 
where $\xi =(t,s)$. If $f$ and $g$ are two $K_{\xi }$-equivalent 
elements in $C^{\xi }_0(M,N)$, that is, there are sequences
$f_n$ and $g_n$ in $C^{\xi }_0(M,N)$ converging to $f$ and $g$ 
respectively and also a sequence $\psi _n \in Diff^{\xi }_0(M)$
such that $f_n(x)=g_n(\psi _n(x))$ for each $x\in M$, then
$\pi _k(f_n)=:f_{n,k}$ and $g_{n,k}:=\pi _k(g_n)$ converge
to $\pi _k(f)$ and $\pi _k(g)$ respectively and also
$\psi _{n,k}:=\pi _k(\psi _n)\in Hom_0(M_k)$.
From the equality $f_{n,k}(x(k))=g_{n,k}(\psi _{n,k}(x(k)))$
for each $n\in \bf N$ and $x(k)\in M_k$ it follows, that
the equivalence relation $K_{\xi }$ induces the corresponding 
equivalence relation $K_{\xi ,k}$ in $\pi _k(C^t_0(M,N))$
such that classes $<\pi _k(f)>_{K,\xi ,k}$ of 
$K_{\xi ,k}$-equivalent elements are closed.
Each element $f_k\in \pi _k(C^{\xi }_0(M,N))$
is characterized by the equality $f_k(s_{0,k})=y_{0,k}$.
Each $\Omega (M_k,N_k)$ is the finite discrete set, since
each $N_k^{M_k}$ is the finite discrete set. 
This induces the quotient mapping $\pi _k: \Omega _t(M,N)\to
\Omega (M_k,N_k)$ and surjective mappings $\pi ^l_k: \Omega (M_l,N_l)
\to \Omega (M_k,N_k)$ for each $l\ge k$ that produces the inverse 
sequence of finite discerete spaces, hence the limit of the
inverse sequence is compact and totally disconnected.
It remains to verify that $\Omega ^w(M,N)$  
is the commutative topoogical monoid with the unit
element and the cancellation property.
\par From the equality $M={\bar M}\setminus \{ s_0 \} $, it follows
that $M_k={\bar M}_k,$ since for each $k\in \bf N$ there exists
$x\in M$ such that $x+B({\bf K^{\sf m}},0,p^{-k})\ni s_0$.
Moreover, $M_k$ and $N_k$ are finite dicrete spaces.
Then $\pi _k(M\vee M)=M_k\vee M_k$, where $A\vee B:=
A\times \{ b_0 \} \cup \{ a_0 \} \times B\subset A\times B$
is the wedge product of pointed spaces $(A,a_0)$ and $(B,b_0)$,
$A$ and $B$ are sets with marked points $a_0\in A$ and $b_0\in B$.
The composition operation is defined on threads 
$ \{ <f_k>_{K,\xi ,k}: k \in {\bf N} \} $ of the inverse sequence
in the following way. There was fixed a $C^{\infty }$-diffeomorphism
$\chi : M\vee M\to M$. Let $x\in M$, then $\pi _k(x)\in M_k$
and $\chi ^{-1}(U)\in M\vee M$, where $U:=\pi _k^{-1}(x+B({\bf K},0,
p^{-k})\cap M.$ On the other hand $\chi ^{-1}(U)$ is a disjoint union 
of balls of radius $p^{-2k}$ in $B({\bf K^{2m}},0,1)$, hence
there is defined a surjective mapping $\chi _k: M_{2k}\vee M_{2k}
\to M_k$ induced by $\chi $, $\pi _k$ and $\pi _{2k}$
such that $\chi _k(\chi ^{-1}(U))=\pi _k(x)$.
If $f$ and $g\in C^{\xi }(M,N)$, then $f\vee g\in C^{\xi }((M\vee 
M),N)$ and $\chi (f\vee g)\in C^{\xi }(M,N)$ as in \S 2.6 \cite{lubp}.
Hence $\chi _k(f_{2k}\vee g_{2k})\in C^{\xi }(M_k,N_k)$
and inevitably $\chi _k(<f_{2k}\vee g_{2k}>_{K,\xi ,2k})
=\chi _k(<f_{2k}>_{K,\xi ,2k}\vee <g_{2k}>_{K,\xi ,2k})
\in \Omega (M_k,N_k)$. 
\par There exists a one to one correspondence between
elements $f\in C_w({\bar M},N)$ and 
$ \{ f_k: k \} \in \{ N_k^{M_k}: k \} $.
Therefore, $pr-\lim_k \Omega (M_k,N_k)$
algebraically is the commutative monoid with the cancellation property.
Let $U$ be a neighbourhood of $e$ in $\Omega ^w(M,N)$, then
there exists $U_k=\pi _k^{-1}(V_k)$ such that $V_k$ is open in
$\Omega (M_k,N_k)$, $e\in U_k$ and $U_k\subset U$. 
On the other hand there exists $U_{2k}=\pi _{2k}^{-1}(V_{2k})$
such that $V_{2k}$ is open in $\Omega (M_{2k},N_{2k})$,
$e\in U_{2k}$ and $U_{2k}+U_{2k}\subset U_k$.
Therefore, $(f+U_{2k})+(g+U_{2k})\subset f+g+U_k\subset f+g+U$
for each $f, g\in \Omega ^w(M,N)$, consequently,
the composition in $\Omega ^w(M,N)$ is continuous.
Since $C^{\xi }_0(M,N)$ is dense in $C_{0,w}({\bar M},N)$,
then $\Omega _{\xi }(M,N)$ is dense in $\Omega ^w(M,N)$.
\par {\bf 3.2. Note.} The compactificiation 
of $\Omega _{\xi }(M,N)$ given above
is not unique. Another compactification is given below.
The second is larger than the first one.
Using the Grothendieck construction we get a compactification
$L^w(M,N)={\bar F}/{\bar B}$ of a loop group $L_{\xi }(M,N)$, where
$\bar F$ is a closure in $(\Omega ^w(M,N))^{\bf Z}$
of a free commutative group $F$ generated by $\Omega ^w(M,N)$ 
and $\bar B$ is a closure of a subgroup $B$ generated by 
all elements $[a+b]-[a]-[b]$, since the product of compact spaces 
is compact by the Tychonoff theorem.
\par {\bf 3.3.} Let now $s_0=0$ and $y_0=0$ be two marked points in 
the compact manifolds $\bar M$ and $N$ embedded into
$\bf K^{\sf m}$ and $\bf K^{\sf n}$ respectively.
There is defined the following $C^{\infty }$-diffeomorphism $inv: 
({\bf K^{\sf m}})' \to ({\bf K^{\sf m}})'$ for
$({\bf K^{\sf m}})':={\bf K^{\sf m}}\setminus 
\{ x:$ $\mbox{ there exists }$ $j
\mbox{ with }$ $x_j=0 \} $ such that $inv (x_1,...,x_{\sf m})=
(x_1^{-1},...,x_{\sf m}^{-1})$. 
Let $M'=M\cap ({\bf K^{\sf m}})'$, then $inv (M')$ is locally compact and 
unbounded in $\bf K^{\sf m}$, consequently, $\pi _k (inv (M'))=(inv (M'))_k$
is a discrete infinite subset in $\bf K_k^{\sf m}$ for each $k\in \bf N$.  
Analogously $\pi _k (inv (M'\vee M'))=(inv (M'\vee M'))_k
\subset \bf K_k^{\sf 2m}$. There exists a $C^{\infty }$-diffeomorphism
$\chi : M\vee M\to M$ such that $inv\circ \chi \circ inv$
is the $C^{\infty }$-diffeomorphism of $inv (M'\vee M')$ with $inv (M')$
and it induces bijective mappings $\chi _k$ of $inv ((inv (M'\vee M'))_k)$
with $inv ((inv (M'))_k)$ for each $k\in \bf N$ such that 
${\hat \pi }^l_k\circ \chi _l=\chi _k$ for each $l\ge k$, where
${\hat \pi }^l_k:=inv\circ \pi ^l_k\circ inv $. 
This produces inverse sequences of discrete spaces
$inv ((inv (M'))_k)=:{\hat M}_k$, $inv ((inv (M'\vee M'))_k)=
{\hat M}_k\vee {\hat M}_k$ and their bijections $\chi _k$
such that $pr-\lim_k{\hat M}_k$ is homeomorphic with $M'$
and $pr-\lim_k \chi _k$ is equal to $\chi $ up to the homeomorphism,
since $pr-\lim_k{\bf K_k^{\sf m}}=\bf K^{\sf m}$ (see also
about admissible modifications and polyhedral expansions
in \cite{luumpe}). If $\psi \in Diff^{\xi }_0({\bar M})$, then
${\hat \psi }\in Diff^{\xi }({\hat M})$. 
Let $J_{f,k}:= \{ h_k: h_k=f_k\circ \psi _k,
\psi _k\in Hom({\hat M}_k), \psi _k(s_{0,k})=s_{0,k} \} $ for
$f_k\in N_k^{{\hat M}_k}$ with $\lim_{x\to 0}f_k(x)=0$,
then $J_{f,k}$ is closed and ${\hat \pi }_k(<f>_{K,\xi })\subset J_{f,k}$.
Therefore, $g_k$ and $f_k$ are ${\hat K}_{\xi ,k}$-equivalent
if and only if there exists $\psi _k\in Hom({\hat M}_k)$ 
such that $\psi _k(s_{0,k})=s_{0,k}$ and $g_k(x)=f_k(\psi _k(x))$ 
for each $x\in {\hat M}_k$. 
Let $\Omega ({\hat M}_k,N_k):={\hat \pi }_k(\Omega _{\xi }(M,N))$.
\par {\bf Theorem.} {\it The set of $\Omega ({\hat M}_k,N_k)$
forms an inverse sequence $S =\{ \Omega ({\hat M}_k,N_k); 
{\hat \pi }^l_k; k\in {\bf N} \} $ such that $pr-\lim S=:
\Omega ^{i,w}(M,N)$ is an associative topological loop monoid with the 
cancellation property and the unit element $e$.
There exists an embedding of $\Omega _{\xi }(M,N)$ into
$\Omega ^{i,w}(M,N)$ such that $\Omega _{\xi }(M,N)$
is dense in $\Omega ^{i,w}(M,N)$.}
\par {\bf Proof.} Let ${U'}_i$ be an analytic disjoint atlas
of $inv(M')$, $f\in C^{\xi }(inv(M'),{\bf K})$, $\psi 
\in Diff^{\xi }(inv(M'))$, then each restriction 
$f|_{{U'}_i}$ has the form $f|_{{U'}_i}(x)=\sum_mf_{i,m}{\bar Q}_{i,m}(x)$
for each $x\in {U'}_i$, where ${\bar Q}_{i,m}$ are basic Amice polynomials
for ${U'}_i$, $f_{i,m}\in \bf K$. Therefore $f$ is a combination
$f=\nabla _if|_{{U'}_i}$, hence ${\hat \pi }_k(f\circ \psi (x))=
\sum_m[{\hat \pi }_k(f_{i,m})\nabla _{(i,\psi _k(x(k))
\in {\hat \pi }_k({U'}_k)}{\bar Q}_{i,m,k}(\psi _k(x(k)))]$ and inevitably
${\hat \pi }_k((f\circ \psi )(x))=f_k\circ \psi _k(x(k))$, where
${\bar Q}_{i,m,k}:={\hat \pi }_k({\bar Q}_{i,m})$, $x\in inv(M')$ 
and $x(k)={\hat \pi }_k(x)$. 
\par As in \S 2.6.2 \cite{lubp} we choose an infinite atlas
$At'(M):= \{ ({U'}_j,{\phi '}_j): j\in {\bf N} \} $ such that
${\phi '}_j: {U'}_j\to B(X,{y'}_j,{r'}_j)$ are homeomorphisms,
$\lim_{k\to \infty }{r'}_{j(k)}=0$, $\lim_{k\to \infty }
{y'}_{j(k)}=0$ for an infinite sequence $\{ j(k)\in {\bf N}: 
k\in {\bf N} \} $ such that 
$cl_{\bar M}[\bigcup_{k=1}^{\infty }{U'}_{j(k)}]$ is a clopen neighbourhood 
of zero in $\bar M$, where $cl_{\bar M}A$ denotes the closure of a subset
$A$ in $\bar M$. We take $|{y'}_{j(k)}|>{r'}_{j(k)}$ for each $k$,
hence $inv (B(X,{y'}_j,{r'}_j)\cap X')=B(X,{y'}_j^{-1},{r'}_j^{-1})\cap X'$
and $\bigcup_k inv({U'}_{j(k)}\cap X')$ is open in $X'$,
where $X={\bf K}^{\sf m}$.
For an atlas $At'(M\vee M):=\{ (W_l,\xi _l): l\in {\bf N} \} $ 
with homeomorphisms $\xi _l: W_l\to B(X,z_l,a_l)$, $\lim_{k\to \infty }
a_{l(k)}=0$, $\lim_{k\to \infty }z_{l(k)}=0$ for an infinite sequence
$\{ l(k)\in {\bf N}: k\in {\bf N} \} $ such that
$cl_{\bar M\vee \bar M}[\bigcup_{k=1}^{\infty }W_{l(k)}]$ 
is a clopen neighbourhood of $0\times 0$ in $\bar M\vee \bar M$
we also choose $|z_l|>a_l$ for each $l$, where
$card ( {\bf N}\setminus \{ l(k): k\in {\bf N} \} )
=card ( {\bf N}\setminus \{ j(k): k\in {\bf N} \} )$.
Then we take $\chi (W_{l(k)})={U'}_{j(k)}$ for each $k\in \bf N$
and $\chi (W_l)={U'}_{\kappa (l)}$ for each 
$l\in ({\bf N}\setminus \{ l(k): k\in {\bf N} \} )$,
where $\kappa : ({\bf N}\setminus \{ l(k): k\in {\bf N} \} )
\to ({\bf N}\setminus \{ j(k): k\in {\bf N} \} )$
is a bijective mapping such that
$p^{-1}\le {r'}_{j(k)}/a_{l(k)}\le p$ for each $k$ and
$p^{-1}\le {r'}_{\kappa (l)}/a_l\le p$ for each
$l\in ({\bf N}\setminus \{ l(k): k\in {\bf N} \} )$.
We can choose the locally affine mapping $\chi $ such that
$\Phi ^n\chi =0$ for each $n\ge 2$ and $B(X',{y'}_l^{-1},{r'}_l^{-1})$
are diffeomorphic with $inv ({U'}_l\cap X')$ and 
$B(X'\vee X',z_l^{-1},a_l^{-1})$ are diffeomorphic with
$inv (W_l\cap (X'\vee X'))$. 
\par This induces the diffeomorphisms ${\hat \chi }:=
inv\circ \chi \circ inv: {\hat M}\vee 
{\hat M}\to {\hat M}$ and ${\hat \chi }^*: C^{\xi }_0(({\hat M}
\vee {\hat M}, \infty \times \infty ), (N,y_0))\to
C^{\xi }_0(({\hat M},\infty ),(N,y_0))$,
since each $\Phi ^n(f\vee g)({\hat \chi }^{-1})$ has an expression
through $\Phi ^l(f\vee g)$ and $\Phi ^j({\hat \chi }^{-1})$
with $l, j\le n$ and $n$ subordinated to $\xi $,
where ${\hat M}:= inv (M')$ and conditions defining the subspace
$C^{\xi }_0 (({\hat M},\infty ),(N,y_0))$ differ from that of
$C^{\xi }_0 ((M,s_0),(N,y_0))$ by substitution of
$\lim_{x\to s_0}$ on $\lim_{|x|\to \infty }$.
Then $\lim_{|x|\to \infty }|{\hat \chi }(x)|=\infty $,
consequenlty, there exists $k_0\in \bf N$
such that ${\hat \chi }_k: {\hat M}_k\vee {\hat M}_k
\to {\hat M}_k$ are bijections for each $k\ge k_0$,
where ${\hat \chi }_k:={\hat \pi }_k\circ {\hat \chi }$.
If $\psi \in Diff^{\xi }(\bar M)$ and $\psi (0)=0$, then
$\lim_{|x|\to \infty }{\hat \psi }(x)=\infty $
and $\lim_{|x|\to \infty }{\hat \psi }^{-1}(x)=\infty .$
Then considering ${\hat \psi }_k$ we get an equivalence
relation $K_{\xi ,k}$ in 
$\{ f_k: f_k\in N_k^{{\hat M}_k}, \lim_{|x|\to \infty }
f_k(x)=0 \} $ induced by $K_{\xi }$, 
where ${\hat M}_k$ is supplied with the quotient norm
induced from the space $X$, since $X'\subset X$, $x\in {\hat M}_k$.
Let $J_k$ denotes the quotient mapping corresponding to $K_{\xi ,k}.$
Therefore analogously to \S 2.6 \cite{lubp}
we get, that $\Omega ({\hat M}_k,N_k)$ are commuative monoids 
with the cancellation property and the unit elements $e_k$, since
$\Omega ({\hat M}_k,N_k)=
\{ f_k: f_k\in C^0({\hat M}_k,N_k),
\lim_{|x|\to \infty }f_k(x)=0 \} /{\hat K}_{\xi ,k}$ and mappings 
${\hat \pi }^l_k: ({{\bf K^{\sf m}})'}_l\to ({{\bf K^{\sf m}})'}_k$
and mappings $\pi ^l_k: {\bf K^{\sf n}}_l\to {\bf K^{\sf n}}_k$
induce mappings ${\hat \pi }^l_k: \Omega ({\hat M}_l,N_l)
\to \Omega ({\hat M}_k,N_k)$ for each $l\ge k$.
Let the topology in $\{ f_k: f_k\in C^0({\hat M}_k,N_k),
\lim_{|x|\to \infty }f_k(x)=0 \} $ be induced from the Tychonoff 
product topology in $N_k^{{\hat M}_k}$ and 
$\Omega ({\hat M}_k,N_k)$ be in the quotient topology.
The space $N_k^{{\hat M}_k}$ is metrizable by the Baire metric
$\rho (x,y):=p^{-j}$, where $j=\min \{ i: x_i\ne y_i,
x_1=y_1,...,x_{i-1}=y_{i-1} \} $, $x=(x_l: x_l\in N_k, l\in {\bf N} )$,
${\hat M}_k$ as enumerated as $\bf N$.
Therefore, $\Omega ({\hat M}_k,N_k)$ is metrizable and
the mapping $(f_k,g_k)\to f_k\vee g_k$ is continuous,
hence the mapping $(J_k(f_k),J_k(g_k))\to J_k(f_k)\circ J_k(g_k)$
is also continuous. Then $J_k(w_{0,k})$ is the unit element,
where $w_{0,k}({\hat M}_k)=0$.
Hence $\Omega ^{i,w}(M,N)$ is the commutative monoid with the cancellation
property and the unit element. Certainly $\prod_k\Omega ({\hat M}_k,
N_k)$ is the topological monoid and $pr-\lim S$ is a closed in it
topological totally disconnected monoid.
For each $f\in C^{\xi }_0(M,N)$ there exists an inverse sequence
$\{ f_k: f_k={\hat \pi }_k(f), k\in {\bf N} \} $
such that $f(x)=pr-\lim_kf_k(x(k))$ for each $x\in M'$, but $M'$ 
is dense in $M$. Therefore there exists an embedding
$\Omega ^{\xi }(M,N)\hookrightarrow \Omega ^{i,w}(M,N)$.
Since $C^{\xi }(M,N)$ is dense in $C^0_0(M,N)$, then  
$\Omega ^{\xi }(M,N)$ is dense in $\Omega ^{i,w}(M,N)$.
\par {\bf 3.4. Corollary.} {\it The inverse sequence of loop monoids induces
the inverse sequence of loop groups $S_L:=\{ L({\hat M}_k,N_k);
{\hat \pi }^l_k; {\bf N} \} $. Its projective limit $L^{i,w}(M,N):=
pr-\lim S_L$ is a commutative topological 
totally disconnected group and $L_{\xi }(M,N)$ has
an embedding in it as a dense subgroup.}
\par {\bf Proof.} Due to Grothendieck construction the inversion operation
$f_k\mapsto f_k^{-1}$ is continous in $L({\hat M}_k,N_k)$
and homomorphisms ${\hat \pi }^l_k$ and ${\hat \pi }_k$ 
have continous extensions
from loop submonoids onto loop groups $L({\hat M}_k,N_k)$.
Each monoid $\Omega ({\hat M}_k,N_k)$ is totally disconnected,
since $N_k^{{\hat M}_k}$ is totally disconnected and
$\Omega ({\hat M}_k,N_k)$ is supplied with the quotient ultrametric,
hence the free Abelian group $F_k$ generated by $\Omega ({\hat M}_k,N_k)$
is also totally disconnected and ultramertizable, consequenlty,
$L({\hat M}_k,N_k)$ is ultrametrizable. Evidently their inverse
limit is also ultrametrizable and the equivalent ultrametric can be chosen 
with values in ${\tilde \Gamma }_{\bf K}:=\{ |z|: z\in {\bf K} \} $,
where ${\tilde \Gamma }_{\bf K}\cap (0,\infty )$ is discrete in
$(0,\infty ):=\{ x: 0<x<\infty , x\in {\bf R} \} $.
Then the projective limit (that is, weak) topology of $L^{i,w}(M,N)$
is induced by the weak topology of $C^0(M,{\bf K})$.
When $M$ and $N$ are non-trivial, then certainly this weak topology 
is strictly weaker, than that of $L_0(M,N)$.
\par {\bf 3.5 Theorem.} {\it For each prime number $p$
the loop group $L_{\xi }(M,N)$ in its weak topology inherited from
$L^{i,w}(M,N)$ has the $p$-adic completion
isomorphic with ${\bf Z_p}^{\aleph _0}.$}
\par {\bf Proof.} If $\bf P$ is an extension of 
the field $\bf K$, then the projective ring homomorphism 
$\pi _k: {\bf K}\to {\bf K_k}$
induces the ring homomorphism $\pi _k : {\bf P}\to 
{\bf P_k}$, then ${\hat \pi }_k ( {\bar \Phi }^v(f(x;h_1,...,h_n;
\zeta _1,...,\zeta _n))={\bar \Phi }^vf_k(x(k);h_1(k),...,h_n(k);
\zeta _1(k),...,\zeta _n(k))$, where $n=[v]+ sign \{ v \} $,
$j_b: {\bf K}\to {\bf P}$, ${\bar \Phi }^vf_k$ is defined for 
the field of fractions generated by ${\bf P_k}$ (see also
\S \S 2.1-2.6 \cite{lubp}). Then the condition
$$\lim_{|x|\to \infty }{\bar \Phi }^vf(x;h_1,...,h_n;
\zeta _1,...,\zeta _n)=0$$ implies the condition
$$\lim_{|x(k)|\to \infty }{\bar \Phi }^vf_k(x(k);h_1(k),...,h_n(k);
\zeta _1(k),...,\zeta _n(k))=0.$$
Therefore, $supp (f_k):={\hat M}_k^f:= \{ x(k): f_k(x(k))\ne 0 \} $ 
is a finite subset of the discrete space ${\hat M}_k$
for each $k\in \bf N$. Then evidently,
${\hat \pi }_k(<g>_{K,\xi })$ is a closed subset in $N_k^{{\hat M}_k}$
for each $g\in C^{\xi }_0(({\hat M},\infty ),(N,0))$,
since for each limit point $f_k$ of ${\hat \pi }_k(<g>_{K,\xi })$ 
its support is the finite subset in ${\hat M}_k$.
Let $k_0$ be such that $N_{k_0}\ne \{ 0 \} $, then this is also true
for each $k\ge k_0$. If $f_k\notin {\hat \pi }_k(<w_0>_{K,\xi })$
and $k\ge k_0$, then $f_k^{\vee n}\notin {\hat \pi }_k(<w_0>_{K,\xi })$
for each $n\in \bf N$, where $f_k^{\vee n}:=f_k\vee ... \vee f_k$
denotes the $n$-times wedge product, since 
$\| f^{\vee n}\|_{C^{\xi }} \ge \| f \|_{C^{\xi }} >0$ and  
$\| f_k^{\vee n}\|_{C({\bf K^{\sf m}_k},{\bf K^{\sf n}})} 
\ge \| f \|_{C({\bf K^{\sf m}_k},{\bf K^{\sf n}})} >0,$ 
where $C({\bf K^{\sf m}_k},{\bf K^{\sf n}})=\pi _k
(C^{\xi }({\bf K^{\sf m}},{\bf K^{\sf n}}))$ is the quotient 
module over the ring $\bf K_k$. Each ${\hat \pi }_k(<f>_{K,\xi })$ 
can be presented as the following composition $v_1b_1+...+v_lb_l$ 
in the additive 
group $L({\hat M}_k,N_k)$, where each $b_i$ corresponds to 
${\hat \pi }_k(<g_i>_{K,\xi })$
and the embedding of $\Omega ({\hat M}_k,N_k)$ into $L({\hat M}_k,N_k)$, 
$v_i\in \{ -1, 0, 1 \} $,
$l=card ({\hat M}_k^f)$, ${\hat M}_k^{g_i}$ are singletons for each
$i=1,...,l$. Using the group $Hom_0(N_k)$ we get that
$L({\hat M}_k,N_k)$ is isomorphic with ${\bf Z}^{n_k}$, where
$n_k=card(N_k)>1$. In view of Corollary 3.4 $L_{\xi }(M,N)$ has
the $p$-adic completion isomorphic with ${\bf Z}_p^{\aleph _0}$,
since $\bf Z$ is dense in $\bf Z_p$ and $pr-\lim_k{\bf Z}^{n_k}=
{\bf Z}^{\aleph _0}.$
\par {\bf 3.6. Note.} 
Using quotient mappings $\eta _{p,s}: {\bf Z}\to {\bf Z}/p^s\bf Z$ 
we get that $L_{\xi }(M,N)^{\aleph _0}$ has the compactification 
equal to $\prod_{p\in {\sf P}}{\bf Z_p}^{\aleph _0},$
where $\sf P$ denotes the set of all prime numbers $p>1$, $s\in \bf N$.
These compactifications produce characters of $L_{\xi }(M,N)$,
since each compact Abelian group has only one-dimensional irreducible 
unitary representations \cite{hew}. On the other hand, there are
irreducible continuous representations of compact groups
in non-Archimedean Banach spaces \cite{roosch}. Among them there are 
infinite-dimensional \cite{diar}. Moreover, in their initial 
topologies diffeomorphism and loop groups also have infinite-dimensional 
irreducible unitary representations \cite{lutmf,lubp}.
\par The problem about $p$-adic completions of diffeomorphism 
and loop groups of manifolds on non-Archimedean Banach spaces
over local fields was formulated by B. Diarra 
after reading articles of S.V. Ludkovsky
on such groups. Then S.V. Ludkovsky investigated this problem and 
all his results and proofs were thoroughly and helpfully discussed
and were corrected with B. Diarra.
 
\par
\par Addresses: S.V. Ludkovsky, Theoretical Department,
\par Institute of General Physics,
\par Russian Academy of Sciences,
\par Str. Vavilov 38, Moscow, 117607, Russia; \\
\par B. Diarra, D\'epartement de Math\'ematiques,
\par Complexe Scientifique des C\'ezeaux,
\par 63177 Aubi\`ere CEDEX, France.

\begin{thebibliography}{99}
\bibitem{ami} Y. Amice. "Interpolation $p$-Adique". Bull. Soc. Math. France
{\bf 92}(1964), 117-180.
\bibitem{aref} Aref'eva I.Ya., Dragovich B., Frampton P.H.,
Volovich I.V. "Wave functions of the universe and $p$-adic gravity".
Int. J. Modern Phys. {\bf 6} (1991), 4341-4358.
\bibitem{arsch} J. Araujo, W.H. Schikhof. "The Weierstrass-Stone
approximation theorem for $p$-adic $C^n$-functions". Ann. Math.
Blaise Pascal. {\bf 1} (1994), 61-74.
\bibitem{boum} N. Bourbaki. "Vari\'et\'es diff\'erentielles et analytiques".
Fasc. XXXIII (Paris: Hermann, 1967).
\bibitem{diar} B. Diarra. "On reducibility of ultrametric almost periodic
linear representations". Glasgow Math. J. {\bf 37} (1995), 83-98.
\bibitem{eng} R. Engelking. "General topology". Second Edit., 
Sigma Ser. in Pure Math. V. {\bf 6} (Berlin: Heldermann Verlag, 1989).
\bibitem{hew} E. Hewitt, K.A. Ross. "Abstract harmonic analysis"
(Berlin: Springer, 1979).
\bibitem{litl} D.E. Littlewood. "The theory of group characters
and matrix representations of groups" (Oxford, 1950).
\bibitem{luseabm} S.V. Ludkovsky. "Irreducible Unitary 
Representations of Non-Archimedean Groups of Diffeomorphisms".
Southeast Asian Bulletin of Mathematics. {\bf 22} (1998),
419-436.
\bibitem{lutmf} S.V. Ludkovsky. "Measures on Diffeomorphism
Groups for Non-Archimedean Manifolds: Group Representations 
and Their Applications". Theoret. and Math. Phys.
{\bf 119: 3} (1999), 698-711.
\bibitem{lubp} S.V. Ludkovsky. 
"Quasi-invariant measures on non-Archimedean
groups and semigroups of loops and paths, their representations".
Annales Mathematiques Blaise Pascal.
{\bf 7: 2} (2000), 19-53, 55-80.
\bibitem{luumpe} S.V. Ludkovsky. "Non-Archimedean Polyhedral Expansions 
of Ultrauniform Spaces".
Russ. Math. Surveys. V. {\bf 54: 5} (1999), 163-164. 
(in detail: Los Alamos National Laboratory, USA.
Preprint {\bf math.AT/0005205}, 39 pages, 22 May 2000).
\bibitem{luum985} S.V. Ludkovsky. "Embeddings of Non-Archimedean Banach 
Manifolds into Non-Archimedean Banach Spaces".
Russ. Math. Surv. {\bf 53} (1998), 1097-1098.
\bibitem{luladg} S.V. Ludkovsky.
"A Structure and Representations of Diffeomorphism Groups
of Non-Archimedean Manifolds". Los Alamos National Laboratory, USA.
Preprint {\bf math.GR/0004126}, 32 pages, 19 April 2000.
\bibitem{mensk} M.B.  Mensky.  "The Path Group.
Measurements.  Fields.  Particles"  (Moscow: Nauka, 1983).  
\bibitem{nari} L. Narici, E. Beckenstein. "Topological vector spaces"
(New York: Marcel Dekker Inc., 1985).
\bibitem{roo} A.C.M. van Rooij. "Non-Archimedean Functional Analysis"
(New York: Marcel Dekker Inc., 1978).
\bibitem{roosch} A.C.M. van Rooij, W.H. Schikhof.
"Groups representations in non-Archimedean Banach spaces".
Bull. Soc. Math. France. Memoire. {\bf 39-40} (1974), 329-340.
\bibitem{sch1} W.H. Schikhof. "Ultrametric calculus"
(Cambridge: Camb. Univ. Press, 1984).
\bibitem{wei} A. Weil. "Basic number theory" (Berlin:Springer-Verlag,1973).
\bibitem{weyl} H. Weyl. "Classical groups, their invariants and 
representations" (Moscow, Inostr. Lit., 1947).
\end{thebibliography}
\end{document}